\documentclass[EJP]{ejpecp} 

\newcommand{\R}{\mathbb{R}}

\newcommand{\N}{\mathbb{N}}
\newcommand{\Q}{\mathbb{Q}}
\newcommand{\E}{\mathds{E}}
\newcommand{\Pb}{\mathds{P}}
\newcommand{\ind}{\mathds{1}}

\def\S{\mathbb{S}}
\def\T{\mathbb{T}}

\def\md{\mid}
\def\Bb#1#2{{\def\md{\bigm| }#1\bigl[#2\bigr]}}

\def\Eb{\Bb\E}

\def\bi{\begin{itemize}}
\def\ei{\end{itemize}}

\def\bnum{\begin{enumerate}}
\def\enum{\end{enumerate}}
\def\LB{\mathcal{B}} 
\def\<#1{\langle #1 \rangle}

\def\p{\mathbf{p}}
\def\t{\mathbf{t}}
\def\r{\mathbf{r}}

\def\s{\mathbf{s}}

\SHORTTITLE{On the heat kernel and the Dirichlet form of Liouville Brownian motion} 

\TITLE{On the heat kernel and the Dirichlet form of~Liouville~Brownian~motion\thanks{Support: ANR grant MAC2 10-BLAN-0123; ANR grant ANR-11-JCJC CHAMU.}} 

\AUTHORS{%
  Christophe Garban\footnote{\'Ecole Normale Sup\'erieure de Lyon (UMPA) and
    CNRS, 69364 Lyon, France.}
  \and 
  R\'emi Rhodes\footnote{Universit{\'e} Paris-Dauphine, Ceremade, F-75016
    Paris, France.} \and Vincent Vargas\footnotemark[3]}


\AMSSUBJ{60G60 ; 60G15 ; 28A80} 

\SUBMITTED{August 1, 2013} 
\ACCEPTED{May 29, 2014} 

\VOLUME{19}
\YEAR{2014}
\PAPERNUM{96}
\DOI{v19-2950}

\ABSTRACT{In \cite{GRV}, a Feller process called {\it Liouville Brownian motion} on $\R^2$ has been introduced. It can be seen as a Brownian motion evolving in a random geometry given formally by the exponential of a (massive) Gaussian Free Field $e^{\gamma\, X}$ and is the right diffusion process to consider regarding $2d$-Liouville quantum gravity.  In this note, we discuss the construction of the associated  Dirichlet form, following essentially \cite{fuku} and the techniques introduced in \cite{GRV}. Then we carry  out the analysis of the Liouville resolvent. In particular, we prove that it is strong Feller, thus obtaining the existence of  the {\it Liouville heat kernel} via a non-trivial theorem of Fukushima and al.
 
  One of the motivations which led to introduce the Liouville Brownian motion in \cite{GRV} was to investigate the puzzling Liouville metric through the eyes of this new stochastic process. In particular,  the theory developed for example in \cite{stollmann,sturm1,sturm2}, whose aim is to capture the ``geometry'' of the underlying space out of the Dirichlet form of a process living on that space, suggests a notion of distance associated to a Dirichlet form. More precisely, under some mild hypothesis on the regularity of the Dirichlet form, they provide a distance in the wide sense, called {\it intrinsic metric}, which is interpreted as an extension of Riemannian  geometry applicable to non differential structures. We prove  that the needed mild hypotheses are satisfied but that the associated intrinsic metric unfortunately vanishes, thus showing that renormalization theory remains out of reach of 
 the metric aspect  of Dirichlet forms. }

\begin{document}



This paper is concerned with the study of a Feller process, called the {\it Liouville Brownian motion}, that has been introduced in \cite{GRV}   to have further insight into the geometry of $2d$-Liouville quantum gravity (see also \cite{berest} where the author constructs {\it Liouville Brownian motion} starting from one point). More precisely, one major mathematical problem in (critical) $2d$-Liouville quantum gravity is to construct a random metric on a two dimensional Riemannian manifold $D$, say a domain of $\R^2$ (or the sphere) equipped with the Euclidean metric $dz^2$, which takes on the form
\begin{equation}\label{i.metric}
e^{\gamma X(z)}dz^2
\end{equation}
where $X$ is a (massive) Gaussian Free Field (GFF) on the manifold $D$   and $\gamma\in [0,2[$ is a coupling constant  (see \cite{cf:KPZ,cf:Da,DFGZ,DistKa,GM,Nak} for further details and insights in Liouville quantum gravity). If it exists, this metric should generate several geometric objects: instead of listing them all, let us just say that each object that can be associated to a smooth Riemannian geometry raises an equivalent question in $2d$-Liouville quantum gravity. Mathematical difficulties originate from the short scale logarithmically divergent behaviour of the correlation function of the GFF $X$. So, for  each object that one wishes to define, one has to apply a renormalization procedure. 

For instance, one can define the volume form associated to this metric. The theory of renormalization for measures formally corresponding to the exponential of Gaussian fields with logarithmic correlations first appeared in the beautiful paper \cite{cf:Kah} under the name of {\it Gaussian multiplicative chaos} and applies to the Free Fields. Thereafter, convolution techniques were developed in \cite{cf:DuSh,cf:RoVa1,cf:RoVa} (see also \cite{review} for further references). This allows to make sense of measures formally defined by:
\begin{equation}\label{i.measure}
M(A)=\int_A e^{\gamma X(z)-\frac{\gamma^2}{2}\E[X(z)^2]}\,dz,
\end{equation}
where $dz$ stands for the volume form (Lebesgue measure) on $D$. To be exhaustive, in the case of Gaussian Free fields, one should integrate against $h(z)\,dz$ where $h$ is a deterministic function involving the conformal radius at $z$ but, first, this term does not play an important role for our concerns and, second, may be handled as well with Kahane's theory. This approach was used in \cite{cf:DuSh,Rnew10} (see also \cite{Rnew4,Benj,Rnew7,Rnew12}) to formulate a rigorous and measure based interpretation of the  Knizhnik-Polyakov-Zamolodchikov formula (KPZ for short) originally  derived in \cite{cf:KPZ}. 

In \cite{GRV} (see also \cite{berest} for a construction starting from one point), the authors defined the {\it Liouville Brownian motion}. It can be thought of as the diffusion process associated to the metric \eqref{i.metric} and is formally the solution of the stochastic differential equation:
\begin{equation}\label{i.LBM}
\left\lbrace
\begin{array}{l}
\LB^{x}_{\t=0}=x \\
d\LB^{
x}_\t= e^{-\frac {\gamma} 2 X(\LB^{x}_\t)+\frac {\gamma^2} 4 \E[X(\LB^{x}_\t)^2]} \,d\bar B_\t.
\end{array} \right.
\end{equation}
where $\bar B$ is a standard Brownian motion living on $D$. Furthermore, they proved that this Markov process is Feller and generates a strongly continuous semigroup $(P^X_t)_{t\geq 0}$, which is symmetric in $L^2(D,M)$.  In particular, the Liouville Brownian motion preserves the Liouville measure $M$. They also noticed that one can attach to the Liouville semigroup $(P^X_\t)_{\t\geq 0}$  a Dirichlet form by the formula:
\begin{equation}\label{i.dirichlet}
\Sigma(f,f)=\lim_{\t\to 0}\frac{1}{\t}\int_D \big(f(x)-P^X_\t f(x)\big) f(x) M(dx)
\end{equation}
with domain $\mathcal{F}$, which is defined as the set of functions $f\in L^2(\R^2,M)$ for which the above limit exists and is finite. This expression is rather non explicit.
 
The purpose of this paper is to pursue the stochastic analysis of $2d$-Liouville quantum gravity initiated in \cite{GRV}. We denote by $H^1(D,dx)$ the standard Sobolev space:
\begin{equation*}
H^1(D,dx)=\Big\{f\in L^2(D,dx); \nabla f\in L^2(D,dx)\Big\},
\end{equation*}
and by $H^1_{loc}(D,dx)$ the functions which are locally in $H^1(D,dx)$.  First, we will make explicit the Liouville Dirichlet form \eqref{i.dirichlet}, relying on  techniques developed in \cite{fuku,cf:Kah}, more precisely traces of Dirichlet forms and potential theory. Before stating the result, recall that a Dirichlet form $\Sigma$ defined on some domain $\mathcal{F} \subset L^2(D,M)$ is strongly local if for all $u,v \in \mathcal{F}$ with compact support such that $v$ is constant on a neighborhood of the support of $u$ then $\Sigma(u,v)=0$. Recall also that the Dirichlet form $\Sigma$ is regular if, denoting $C_c(D)$ the space of continious functions with compact support, the set $\mathcal{F} \cap C_c(D)$ is dense in  $C_c(D)$ for the uniform norm and dense in $\mathcal{F}$ for the norm induced by the scalar product $f,g \mapsto \Sigma(f,g)+ \int_D f (x)g (x) M(dx)$.

\begin{theorem}\label{i.th.dirform}
For $\gamma\in [0,2[$, the Liouville Dirichlet form $(\Sigma,\mathcal{F})$ takes on the following explicit form: its domain is 
$$\mathcal{F}=\Big\{f\in L^2(D,M)\cap H^1_{loc}(D,dx); \nabla f\in L^2(D,dx)\Big\},$$
 and for all functions $f,g\in\mathcal{F}$:
$$\Sigma(f,g)=\int_{D}\nabla f(x)\cdot \nabla g(x)\,dx.$$
Furthermore, it is strongly local and regular. 
\end{theorem}

Let us stress here that understanding rigorously the above theorem is not obvious since the Liouville measure $M$ and the Lebesgue measure $dx$ are singular. The domain $\mathcal{F}$ is composed of the functions $u\in L^2(D,M)$ such that there exists a function $f \in H^1_{loc}(D,dx)$ satisfying $\nabla f\in L^2(D,dx)$ and $u(x)=f(x)$ for $M(dx)$-almost every $x$. It is a consequence of the general theory developped in \cite{fuku} (see chapter 6) and of the tools developped in  \cite{GRV} that the definition of $(\Sigma,\mathcal{F})$ actually makes sense: indeed, if $f,g$ in $H^1_{loc}(D,dx)$ are such that $f(x)=g(x)$ for $M(dx)$-almost every $x$ then $\nabla f(x)=\nabla g(x)$ for $dx$-almost every $x$.

Then we perform an analysis of the Liouville resolvent family $(R_\lambda^X)_{\lambda > 0}$ defined on the space $C_b(D)$ of bounded continuous functions by:
$$\forall f \in C_b(D),\quad R_\lambda^X f(x)=\int_0^{+\infty}e^{-\lambda t}P^X_tf(x)\,dt.$$
We will prove that this family possesses strong regularizing properties. In particular, if $B_b(D)$ denotes the set of bounded measurable functions,  our two main theorems concerning the resolvent family are:
\begin{theorem}
Almost surely in $X$, for $\gamma\in [0,2[$, the resolvent operator $R_\lambda^X$ is strong Feller in the sense that it maps the set $B_b(D)$ of bounded measurable functions into the set of continuous bounded functions.
\end{theorem}

\begin{theorem}\label{i.holdresolbis}
Assume  $\gamma\in [0,2[$. There is an exponent $\alpha\in(0,1)$ (depending only on $\gamma$), such that, almost surely in $X$,  for all $\lambda>0$ the Liouville resolvent is locally $\alpha$-H\"older. More precisely, for each $R$ and $\lambda_0>0$, we can find a random constant $C_{R,\lambda_0}$ (depending only on the field $X$), which is $\Pb^X$-almost surely finite, such that for all $\lambda\in ]0,\lambda_0]$ and  for all continuous function $f:D\to \R$ vanishing at infinity, $\forall x,y\in B(0,R)$:
$$  |R^X_\lambda f(x)-R^X_\lambda f(y)|\leq \lambda^{-1}C_{R} \|f\|_\infty|x-y|^\alpha.$$
\end{theorem}

As a  consequence, we obtain the existence of the massive Liouville Green functions, which are nothing but the densities of the resolvent operator with respect to the Liouville measure (see Theorem \ref{th:green}).

For symmetric semigroups, Fukushima and al. \cite{fuku} proved the highly non-trivial theorem (see their Theorems 4.1.2 and 4.2.4) which states that absolute continuity of the resolvent family is equivalent to absolute continuity of the semigroup. As such, this allows us to obtain the following theorem on the existence of a heat kernel:
\begin{theorem}{\bf Liouville heat kernel. }\label{i.th:heat}The Liouville semigroup $(P^X_\t)_{\t>0}$ is absolutely continuous with respect to the Liouville measure. There exists a family $(\p^X_\t(\cdot,\cdot))_{\t\geq 0}$, called the Liouville heat kernel, of jointly measurable functions such that:
$$\forall f \in B_b(D),\quad P^X_\t f(x)=\int_{D}f(y)\p^X_\t(x,y)\,M(dy)$$ and such that:\\
1) (positivity) for all $\t>0$ and for all $x\in D$, for $M(dy)$-almost every $y\in D$, $$\p^X_\t(x,y) \geq 0,$$  
2) (symmetry) for all $\t>0$ and for  every $x,y\in D$: $$\p^X_\t(x,y)=\p^X_\t(y,x),$$
3) (semigroup property) for all $\s,\t\geq 0$, for all $x,y \in D$, $$\p_{\t+\s}^X(x,y)=\int_{D}\p^X_\t(x,z)\p^X_{\s}(z,y)\,M(dz).$$  
\end{theorem}

These properties have interesting consequences regarding the stochastic structure of $2d$-Liouville quantum gravity. For instance, the Liouville Brownian motion spends Lebesgue almost all the times in the set of points supporting the Liouville measure, nowadays called the thick points of the field $X$ and first introduced by Kahane in the case of log-correlated Gaussian fields \cite{cf:Kah} like Free Fields (see also \cite{BMI,hmp}). Furthermore, for a given time $t$, the Liouville Brownian motion is almost surely located on the  thick points of $X$. We will also define the {\it Liouville Green function}   to investigate the ergodic properties of the Liouville Brownian motion, which turns out to be irreducible and recurrent.

Finally, let us end this introduction by a discussion on the Liouville Dirichlet form  as well as its possible relevance to the construction of the Liouville distance.  Over the last 20 years, a rich theory has been developed whose aim is to capture the ``geometry'' of the underlying space out of the Dirichlet form of a process living on that space. See for example \cite{stollmann,sturm1,sturm2}. This geometric aspect of Dirichlet forms can be interpreted in a sense as an extension of Riemannian  geometry applicable to non differential structures. 
Among the recent progresses of Dirichlet forms has emerged the notion of {\it intrinsic metric} associated to  a strongly local regular Dirichlet form \cite{BM1,BM2,DA,stollmann,sturm1,sturm2,VSC}. It is natural to wonder if this theory is well suited to this problem of constructing the Liouville distance. More precisely, the intrinsic metric is defined by
\begin{equation}\label{def.metric}
d_X(x,y)=\sup\{f(x)-f(y);f\in\mathcal{F}_{loc}\cap C(D),\Gamma(f,f)\leq M\}.
\end{equation}
where $\mathcal{F}_{loc}$ is the space of functions which are locally in $\mathcal{F}$ and $\Gamma$ is called the energy measure of $f$ (that will be defined in greater detail in section 4). This distance is actually a distance in the wide sense, meaning that it can possibly take values $d_X(x,y)=0$ or $d_X(x,y)=+\infty$ for some $x\not =y$. Let us point out that, when the field $X$ is smooth enough (and therefore not a free field), the distance \eqref{def.metric} coincides with the Riemannian distance generated by the metric tensor $e^{\gamma X(z)}\,dz^2$. Generally speaking, the point is to prove that the topology associated to this distance is Euclidean, in which case $d_X$ is a proper distance and $(D,d_X)$ is a length space (see \cite[Theorem 5.2]{stollmann}). Unfortunately, in the context of $2d$-Liouville quantum gravity, we prove that this intrinsic metric turns out to be $0$.  Anyway, this fact is also   interesting as it sheds some new light on the mechanisms involved in the renormalization of the Liouville distance (if it exists).  
\subsection{Notations}

We stick to the notations of \cite{GRV} (see  the section "Background"), where the basic tools needed to define $2d$-Liouville quantum gravity are described. In particular, a description of the construction of Free Fields and their cutoff regularization are given: throughout the paper, the field $X$ may thus be a Massive Free Field on the whole plane $D=\R^2$ or a Gaussian Free Field on the $2$-dimensional torus $D=\T^2$ or   sphere $D=\S^2$. $(X_n)_{n}$ stands for the cutoff approximation of $X$ defined in \cite{GRV} and $M$ for the Gaussian multiplicative chaos associated to $X$:
$$M(A)=\int_Ae^{\gamma X(x)-\frac{\gamma^2}{2}\E[X(x)^2]}\,dx,$$ where $\gamma\in [0,2)$ and $A$ is a measurable subset of $D$. 

  \textbf{Also, the reader may find a list of notations used throughout  the paper in Section \ref{index}. In the sequel, we will use these notations with no further notice.}

\section{Liouville Dirichlet form}\label{sec.LFD}

The purpose of this section is to give an explicit description of the Liouville Dirichlet form, namely the Dirichlet form of the Liouville Brownian motion, by combining \cite{fuku} and the results in \cite{GRV}. The first part of this section is devoted to recalling a few material about Dirichlet forms in order to facilitate the reading of this paper. Then we identify the Dirichlet form and, finally, we discuss some questions naturally raised by the construction of the Dirichlet form. Among them: "Can we construct the Liouville Brownian motion from the only use of \cite{fuku}?"

\subsection{Background on positive continuous additive functionals and Revuz measures}

In this subsection, to facilitate the reading of our results, we first summarize the content of section 5 in \cite{fuku} applied to the standard Brownian $(\Omega,(B_t^x)_{t \geq 0},(\mathcal{F}_t)_{t \geq 0},(\Pb^B_x)_{x \in D})$ in $D$ which is of course reversible for the canonical volume form $dx$ of $D$. We suppose that the space $\Omega$ is equipped with the standard shifts $(\theta_t)_{t \geq 0}$ on the trajectory. One may then consider the classical notion of capacity associated to the Brownian motion. In this context, we have the following definitions:

\begin{definition}[Capacity and polar set] 
The capacity of an open set $O\subset D$ is defined by
$${\rm Cap}(O)=\inf\{\int_D|f(x)|^2\,dx+\int_D|\nabla f(x)|^2\,dx;f\in H^1(D,dx),\,\,f \geq 1\text{ over }O \}.$$
The capacity of a Borel measurable set $K$ is then defined as: 
\begin{equation*}
{\rm Cap}(K)= \underset{O \text{open}, K \subset O}{\inf}{\rm Cap}(K). 
\end{equation*}
The set $K$ is said polar when ${\rm Cap}(K)=0$.
\end{definition}

\begin{definition}[Revuz measure] 
A Revuz measure $\mu$ is a Radon measure on $D$ which does not charge the polar sets.
\end{definition}

\begin{definition}[PCAF] 
A positive continuous additive functional (PCAF) $(A_t)_{t \geq 0}$ is a $\mathcal{F}_t$-adapted continuous functional with values in $[0,\infty]$ that satisfies for all $\omega \in \Lambda$:
\begin{equation*}
 A_{t+s}(\omega)=A_s(\omega)+ A_t(\theta_s(\omega)), \quad s,t \geq 0
\end{equation*} 
where $\Lambda$ is defined in the following way: there exists a polar set $N$ (for the standard Brownian motion) such that for all $x \in D \setminus N$, $\Pb^B_x(\Lambda)=1$ and $\theta_t(\Lambda) \subset \Lambda$ for all $t \geq 0$. 

\end{definition}

In particular, a PCAF is defined for all starting points $x \in \R^2$ except possibly on a polar set for the standard Brownian motion. One can also work with a PCAF starting from \textbf{all} points, that is when the set $N$ in the above definition can be chosen to be empty. In that case, the PCAF is said {\it in the strict sense}.

Finally, we conclude with the following definition on the support of a PCAF:
\begin{definition}[support of a PCAF] 
Let $(A_t)_{t \geq 0}$ be a PCAF with associated polar set $N$. The support of $(A_t)_{t \geq 0}$ is defined by:
\begin{equation*}
 \tilde{Y}= \Big\{x \in D  \setminus N: \: \Pb^B_x(R=0)=1\Big\},
\end{equation*}
where $R=\inf \lbrace t>0: \: A_t>0 \rbrace$.
\end{definition}

From section 5 in \cite{fuku}, there is a one to one correspondence between Revuz measures and PCAFs. More precisely, a Revuz measure $\mu$ is associated bijectively to a PCAF $(A_t)_{t \geq 0}$ if the following relation is valid for all nonnegative $f,g \in B_b(D)$ and all $t>0$
\begin{equation*}
\int_{D}\E^B_{x}\Big[\int_0^t f(B^x_s)\,dA_s\Big] g(x)\,dx=\int_0^t \left( \int_{D} \left( \int_D p_s(x,y)f(y) \, \mu(dy) \right ) g(x)\,dx   \right)ds
\end{equation*}  
where $p_s(x,y)$ is the heat kernel of the standard Brownian motion $B$. In the next subsection, we will identify the Liouville measure $M$ as the Revuz measure associated to the increasing functional $F$ constructed in \cite{GRV}. Let us first check that the measure $M$ is a Revuz measure, i.e. it does not charge polar sets:

\begin{lemma}\label{capacity}
Almost surely in $X$, the Liouville measure $M$ does not charge the polar sets of the (standard) Brownian motion.
\end{lemma}

\noindent {\it Proof.}  Let $A$ be a bounded polar set and let $R>0$ be such that $A\subset B(0,R)$. From \cite{revuzp} (see also \cite{kaku}), it suffices to prove that the mapping $x\mapsto \int G_R(x,y)M(dy)$ is bounded, where $G_R$ stands for the Green function of the Brownian motion killed upon touching $\partial B(0,R)$.  Recall that the Green function over $B(0,R)$ takes on the form
$$G_R(x,y)=\ln\frac{1}{d(x,y)}+g(x,y)$$ for some bounded function $g$ over $B(0,R)$, where $d$ stands for the usual Riemannian distance on $D$. 
The result thus follows from \cite{GRV} where it is proved that the Liouville measure uniformly integrates the $\ln$ over compact sets.
\qed

\subsection{The Revuz measure associated to Liouville Brownian motion}
In this subsection, we identify the measure $M$ as the Revuz measure associated to the functional $F$ introduced in \cite{GRV}. This functional $F$  is  defined almost surely in $X$ for all $x\in D$ by
$$F(x,t)=\int_0^te^{\gamma X (B_r^x)-\frac{\gamma^2}{2}\E[X^2(B_r^x)]}\,dr,$$ where $B$ is a standard Brownian motion on $D$.
By setting $$\sigma_x=\inf\{s>0; F(x,s)>0\},$$ it is proven in \cite{GRV} that:
\begin{equation}\label{quasiregular}
\text{a.s. in }X,\forall x\in D,\quad \Pb^{B}_x(\sigma_x=0)=1.
\end{equation}
 We claim:
\begin{lemma}\label{revuz}
Almost surely in $X$, $F$ is a PCAF in the strict sense whose support is the whole domain $D$. Also, the Revuz measure of $F$ is the Liouville measure $M$.
\end{lemma}

\noindent {\it Proof.} 
The fact that $F$ is a PCAF in the strict sense whose support is the whole domain $D$ is a direct consequence of \cite{GRV} as summarized in (\ref{quasiregular}).

The Revuz measure $\mu$ associated to $F$ is  the unique measure on $D$ that does not charge polar sets and such that: 
\begin{equation*}
\int_{D}\E^B_{x}\Big[\int_0^t f(B^x_s)\,dA_s\Big] g(x)\,dx=\int_0^t \left( \int_{D} \left( \int_D p_s(x,y)f(y) \, \mu(dy) \right ) g(x)\,dx   \right)ds
\end{equation*}  
 for all continuous nonnegative compactly supported functions $f,g$.
Here $B^x$ stands for the law of a Brownian motion starting from $x$. Recall that $p_s(x,y)$ is the standard heat kernel on $D$. To identify the measure $\mu$ it suffices to compute its values on the set of continuous functions  with compact support.  For such a function, we have:
\begin{align*}
\E^{B}_x\Big[\int_0^tf(B^x_s)\,F(x,ds)\Big] &= \E^{B}_x\Big[\int_0^tf(B^x_s)\,e^{\gamma X(B^x_s)-\frac{\gamma^2}{2}\E[X(B^x_s)^2]}\,ds\Big] \\
&=\int_0^t\int_{D}f(y)p_s(x,y)e^{\gamma X(y)-\frac{\gamma^2}{2}\E[X(y)^2]}\,dy\,dr\\
&=\int_0^t\int_{D}f(y)p_s(x,y)\,M(dy)\,ds.
\end{align*}
Then 
\begin{equation*}
\int_{D}\E^{B}_x\Big[\int_0^tf(B^x_s)\,F(x,ds)\Big]g(x)\,dx= \int_0^t \left( \int_{D} \left( \int_D p_s(x,y)f(y) \, \mu(dy) \right ) g(x)\,dx   \right)ds
\end{equation*}
The proof is complete.\qed

\subsection{Construction of the Liouville Dirichlet form $(\Sigma, \mathcal{F})$}

In this subsection, we want to apply Theorem 6.2.1 in \cite{fuku}. Recall that the Liouville Brownian motion is defined in \cite{GRV} as a continuous Markov process defined for all starting points $x$ by the relation: $$\mathcal{B}^x_\t=B_{\langle \mathcal{B}^x\rangle_\t}^x$$ 
where $\langle \mathcal{B}^x\rangle$ is defined by
$$\langle \mathcal{B}^x\rangle_\t= \inf \lbrace  s>0; \:  F(x,s) > \t \rbrace. $$
We know that $M$ is the Revuz measure associated to $F$. Hence, we can now straightforwardly apply the abstract framework of Theorem 6.2.1 in \cite{fuku} to get the following expression for the Dirichlet form associated to Liouville Brownian motion:

\begin{theorem}
The Liouville Dirichlet form $(\Sigma,\mathcal{F})$ takes on the following explicit form on $L^2(D,M)$:
\begin{equation}\label{formesimple}
\Sigma(f,g)=\int_{D}\nabla f(x)\cdot \nabla g(x)\,dx
\end{equation} 
with domain
$$\mathcal{F}=\Big\{f\in L^2(D,M)\cap H^1_{loc}(D,dx); \nabla f\in L^2(D,dx)\Big\},$$
Furthermore, it is strongly local and regular.
\end{theorem}

In fact, for any PCAF $(A_t)_{t \geq 0}$ associated to Brownian motion, one can define the Dirichlet form associated to the Hunt process $B_{A_t^{-1}}$. In this general case, theorem 6.2.1 in \cite{fuku} gives an expression to the Dirichlet form which is non explicit and involves an abstract projection construction involving the support $\tilde{Y}$ of the PCAF. Nonetheless, there is one case where the Dirichlet form takes on the simple form (\ref{formesimple}): when the support $\tilde{Y}$  is the whole space $D$ (recall that this constitutes a large part of the work \cite{GRV}).

\begin{remark}
This result may appear surprising for non specialists of Dirichlet forms. Let us forget for a while Liouville quantum gravity and assume that the measure $M$ is a smooth measure, meaning that it has a density w.r.t. the Lebesgue measure bounded from above and away from $0$. Then we obviously have $L^2(D,dx)=L^2(D,M)$. In that case, the domain and  expression of the time changed Dirichlet form coincide with those of the Dirichlet form of the standard Brownian motion on $D$. So, a natural question is: "How do we differentiate the Markov process associated to this time changed Dirichlet form from the standard Brownian motion?". The answer is hidden in the fact that a Dirichlet form uniquely determines a Markovian semi-group provided that you fix a reference measure with respect to which you impose the semi-group to be symmetric. In the case of the standard Brownian motion, the reference measure is the Lebesgue measure $dx$ whereas  the reference measure is $M$ in the case of the time changed Brownian motion.
\end{remark}

\subsection{Discussion about the construction of the Dirichlet form and the associated Hunt process}

A natural question regarding the theory of Dirichlet forms is: "Can one construct directly the Liouville Brownian motion via the theory of Dirichlet forms without using the results  in \cite{GRV}?".  Since the Liouville measure is a  Revuz measure, it uniquely defines a PCAF $(A_t)_t$. This PCAF may be used to change the time of a reference Dirichlet form, here that of the standard Brownian motion on $D$. The time changed Dirichlet form constructed in \cite[Theorem 6.2.1]{fuku} corresponds to that of a Hunt process $H_\t=B_{A_\t^{-1}}$ where $B$ is a standard Brownian motion and $A_\t^{-1}$ is the inverse of the PCAF $(A_t)_t$.  Nevertheless, we stress that identifying  this Hunt process explicitly is not obvious without using the tools developed in \cite{GRV}. Moreover, this abstract construction of $H_\t$ rigorously defines a Hunt process living in the space $D \setminus N$ where $N$ is a polar set. To our knowledge, there is no general theory on Dirichlet forms which enables to get rid of this polar set, hence constructing a PCAF in the strict sense and a Hunt process starting from all points of $D$. In conclusion, without using the tools developped in \cite{GRV}, one can construct the Liouville Brownian motion in a non explicit way living in $D \setminus N$ and for starting points in $D \setminus N$ where $N$ is a polar set (depending on the randomness of $X$); in this context, one can not start the process from one given fixed point $x \in D$ or define a Feller process in the strict sense. Even if this was the case, in order to identify the corresponding Dirichlet form by the simple formula (\ref{formesimple}), one must show that the PCAF has full support (which is also part of the work done in \cite{GRV}).

 Let us mention that a measurable Riemannian structure associated to strongly local regular Dirichlet forms is built in \cite{Hino}. In \cite{Kim}, harmonic functions and Harnack inequalities for trace processes (i.e. associated to time changed Dirichlet forms) are studied. In particular, it is proved that harmonic functions for the Liouville Brownian motion are harmonic for the Euclidean Brownian motion and that harmonic functions for the Liouville Brownian motion satisfy scale invariant Harnack inequalities. Actually, there are many powerful tools that can be associated with a Dirichlet forms and listing them exhaustively  is far beyond the scope of this paper.  

\section{Liouville Heat Kernel and Liouville Green Functions}

The Liouville Brownian motion generates a Feller semi-group $(P^X_\t)_\t$, which can be extended to a strongly continuous semigroup  on $L^p(D,M)$  for $1\leq p<+\infty$ and is reversible with respect to the Liouville measure $M$ (see \cite{GRV}). 
Recall that 
\begin{proposition}[\cite{GRV}]\label{cv:semigroup}
For $\gamma<2$, almost surely in $X$, the $n$-regularized semi-group $(P^n_\t)_\t$  converges towards
the Liouville semi-group $(P^X_\t)_\t$ in the sense that for all function $f\in C_b(D)$:
$$\forall x \in D,\quad \lim_{n\to\infty}P_\t^nf(x)=P^X_\t(x).$$
\end{proposition}

The main purpose of this section is to prove the existence (almost surely in $X$) of a heat-kernel $\p_\t(x,y)$ for this Feller semi-group $(P^X_\t)_{\t \geq 0}$.  Our strategy for establishing the existence of the heat-kernel will be first to prove that the resolvent associated to our Liouville Brownian motion is (a.s. in $X$) absolutely continuous w.r.t the Liouville measure $M$: see Theorem \ref{th:green}. In general, the absolute continuity of the resolvent is far from implying the absolute continuity of the semi-group (think for example of the process defined on the circle by $X_t^x=
e^{i(x + t)}$). Nevertheless, as stated in the introduction, the symmetry of the Liouville semi-group w.r.t. the Liouville measure $M$ allows us to apply a deep theorem of Fukushima and al.   \cite{fuku}  to conclude: see Theorem \ref{th:heat}. Finally we deduce some corollaries along this section such as the fact that the Liouville Brownian motion a.s. spends most of his time in the thick points of the field $X$, the construction of the {\it Liouville Green function} or the study of the ergodic properties of the Liouville Brownian motion.

\subsection{Analysis of the Liouville resolvent and existence of the Liouville heat kernel}

One may also consider the resolvent family $(R^X_\lambda)_{\lambda>0}$ associated to the semigroup $(P^X_\t)_\t$. In a standard way, the resolvent operator reads:
\begin{equation}\label{def:resol}
\forall f\in C_b(D),\quad R^X_\lambda f(x)=\int_0^{\infty}e^{-\lambda \t}P_\t^Xf(x)\,d\t.
\end{equation}
Furthermore, the resolvent family $(R^X_\lambda)_{\lambda>0}$ is self-adjoint in $L^2(D,M)$ and extends  to a strongly continuous resolvent family on  the $L^p(D,M)$ spaces for $1\leq p<+\infty$. This results from the properties of the semi-group. From Proposition \ref{cv:semigroup}, it is straightforward to deduce:
\begin{proposition}\label{cv:resolvent}
For $\gamma<2$, almost surely in $X$, the $n$-regularized resolvent family $(R^n_\lambda)_\lambda$  converges towards
the Liouville resolvent $(R^X_\lambda)_\lambda$ in the sense that for all function $f\in C_b(D)$:
$$\forall x \in D,\quad \lim_{n\to\infty}R^n_\lambda f(x)=R^X_\lambda f(x).$$
\end{proposition}
Also, it is possible to get an explicit expression for the resolvent operator:
\begin{proposition}\label{expr:resol1}
 For $\gamma<2$, almost surely in $X$, the resolvent operator takes on the following form for all measurable bounded function $f$ on $D$: 
 $$R_\lambda^Xf(x)=\E^B_x\big[\int_0^{\infty}e^{-\lambda F(x,t)}f(B^x_t)\,F(x,dt)\big].$$
\end{proposition}

\vspace{1mm}
\noindent {\it Proof.} Given a measurable bounded function $f$ on $D$, we have:
\begin{align*}
R_\lambda^Xf(x)&=\int_0^{\infty}e^{-\lambda t}P^X_tf(x)\,dt\\
&=\int_0^{\infty}e^{-\lambda t}\E^B_x[f(\LB^x_t)]\,dt\\
&=\E^B_x\Big[\int_0^{\infty}e^{-\lambda t}f(B^x_{\langle \LB^x\rangle_t})\,dt\Big]\\
&=\E^B_x\Big[\int_0^{\infty}e^{-\lambda F(x,s)}f(B^x_s)\,F(x,ds)\Big],
\end{align*}
which completes the proof.\qed

\vspace{1mm}
\begin{theorem}\label{expr:resol}
 For $\gamma<2$, almost surely in $X$, the resolvent operator $R_\lambda^X$ is strong Feller, i.e. maps the measurable bounded functions into the set of continuous bounded functions.
\end{theorem}

\vspace{1mm}
\noindent {\it Proof.} Let us   consider a bounded measurable function $f$ and let us prove that $R_\lambda^Xf(x)$ is a continuous function of $x$. To this purpose, write for some arbitrary $\epsilon>0$:
\begin{align*}
R_\lambda^Xf(x)&=\E^{B}_x\Big[\int_0^{\infty}e^{-\lambda F(x,s)}f(B^x_s)\,F(x,ds)\Big]\\
&=\E^{B}_x\Big[\int_0^{\epsilon}e^{-\lambda F(x,s)}f(B^x_s)\,F(x,ds)\Big]+\E^{B}_x\Big[\int_\epsilon^{\infty}e^{-\lambda F(x,s)}f(B^x_s)\,F(x,ds)\Big]\\
&\stackrel{def}{=}N_\epsilon(x)+R_\lambda^{X,\epsilon}f(x).
\end{align*}
We are going to prove that the family of functions $(N_\epsilon)_\epsilon$ uniformly converges towards $0$ on compact subsets of $D$ as $\epsilon\to 0$ (obviously, if $D$ is compact, we will prove  uniform convergence on $D$) and that the functions $R_\lambda^{X,\epsilon}f$ are continuous. First we focus on $(N_\epsilon)_\epsilon$ and write the obvious inequality:
\begin{align*}
\sup_{x\in B(0,R)}|N_\epsilon(x)|\leq & \|f\|_{\infty}\sup_{x\in B(0,R)}\E^{B}_x[F(x,\epsilon)].
\end{align*}
From \cite{GRV}, we know that the latter quantity converges to $0$ as $\epsilon$ goes to $0$. Let us now prove the continuity of $R_\lambda^{X,\epsilon}f$. By the Markov property of the Brownian motion, we get:
\begin{align*}
R_\lambda^{X,\epsilon}f(x)&=\E^{B}_x\Big[\int_\epsilon^{\infty}e^{-\lambda F(x,s)}f(B^x_s)\,F(x,ds)\Big]\\
&=\E^{B}_x\Big[e^{-\lambda F(x,\epsilon)}R_\lambda^Xf(B^x_\epsilon)\Big].
\end{align*}
Now we consider two points $x$ and $y$ in $D$ and realize the coupling of $(B^x,F(x,\cdot))$ and $(B^y,F(y,\cdot)) $ explained in \cite{GRV}. Recall that this coupling lemma allows us to construct  a Brownian motion $B^x$ starting from $x$ and a Brownian motion $B^y$  starting from $y$ in such a way that they coincide after some random stopping time $\tau^{x,y}$. Let us denote by $\Pb^B$ the underlying probability measure and $\E^B$ the corresponding expectation. We obtain:
\begin{align*}
|R_\lambda^{X,\epsilon}&f(x)-R_\lambda^{X,\epsilon}f(y)|\\=&\Big|\E^B\Big[e^{-\lambda F(x,\epsilon)}R_\lambda^Xf(B_\epsilon^x)\Big]-\E^B\Big[e^{-\lambda F(y,\epsilon)}R_\lambda^Xf(B_\epsilon^y)\Big]\Big|\\
\leq &\Big|\E^B\Big[e^{-\lambda F(x,\epsilon)}R_\lambda^Xf(B_\epsilon^x)\Big]-\E^B\Big[e^{-\lambda F(x,\epsilon)}R_\lambda^Xf(B_\epsilon^y)\Big]\Big|\\
&+\Big|\E^B\Big[e^{-\lambda F(x,\epsilon)}R_\lambda^Xf(B_\epsilon^y)\Big]-\E^B\Big[e^{-\lambda F(y,\epsilon)}R_\lambda^Xf(B_\epsilon^y)\Big]\Big|\\
\leq & \E^B\big[\big|R_\lambda^Xf(B_\epsilon^x)-R_\lambda^Xf(B_\epsilon^y)\big|\big]+\lambda^{-1}\|f\|_\infty \E^B\big[\big|e^{-\lambda F(x,\epsilon)}-e^{-\lambda F(y,\epsilon)}\big|\big].
\end{align*}
Concerning the first quantity, observe that it is different from $0$ only if the two Brownian motions have not coupled before time $\epsilon$, in which case we use the rough bound $\|R_\lambda^Xf\|_\infty\leq \lambda^{-1}\|f\|_\infty$ to get:
\begin{align}\label{ineg:easy}
\E^B\big[\big|R_\lambda^Xf(B_\epsilon^x)-R_\lambda^Xf(B_\epsilon^y)\big|\big]&=\E^B\big[\big|R_\lambda^Xf(B_\epsilon^x)-R_\lambda^Xf(B_\epsilon^y)\big|;\epsilon\leq \tau^{x,y}\big]\nonumber\\
&\leq 2\lambda^{-1}\|f\|_\infty\Pb(\epsilon\leq \tau^{x,y}).
\end{align}
This latter quantity converges towards $0$ uniformly on compact sets as $|x-y|\to 0$. The second quantity is treated with the same idea:
\begin{align*}
\E^B\big[&\big|e^{-\lambda F(x,\epsilon)}-e^{-\lambda F(y,\epsilon)}\big|\big]\\
=&\E^B\big[\big|e^{-\lambda F(x,\epsilon)}-e^{-\lambda F(y,\epsilon)}\big|;\epsilon\leq\tau^{x,y}\big]+\E^B\big[\big|e^{-\lambda F(x,\epsilon)}-e^{-\lambda F(y,\epsilon)}\big|;\epsilon>\tau^{x,y}\big]\\
\leq &2\Pb(\epsilon\leq \tau^{x,y})+\E^B\big[\big|e^{-\lambda F(x,\tau^{x,y})-\lambda F(x,]\tau^{x,y},\epsilon])}-e^{-\lambda F(y,\tau^{x,y})-\lambda F(y,]\tau^{x,y},\epsilon])}\big|;\epsilon>\tau^{x,y}\big].
\end{align*}
Observe that, on the event $\{\epsilon>\tau^{x,y}\}$, we have $F(x,]\tau^{x,y},\epsilon])=F(y,]\tau^{x,y},\epsilon])$. We deduce:
\begin{align}
\E^B\big[&\big|e^{-\lambda F(x,\epsilon)}-e^{-\lambda F(y,\epsilon)}\big|\big]\nonumber\\
\leq &2\Pb(\epsilon\leq \tau^{x,y})+\E^B\big[\big|e^{-\lambda F(x,\tau^{x,y})}-e^{-\lambda F(y,\tau^{x,y})}\big|;\epsilon>\tau^{x,y}\big]\nonumber\\
\leq &2\Pb(\epsilon\leq \tau^{x,y})+\E^B\big[\min\big(2,\lambda | F(x,\tau^{x,y})-F(y,\tau^{x,y})|\big)\big]\nonumber\\
\leq &2\Pb(\epsilon\leq \tau^{x,y})+\E^B\big[\min\big(2,\lambda  F(x,\delta)+\lambda F(y,\delta)\big)\big]+2\Pb(\tau^{x,y}>\delta)\label{inegdelta}
\end{align}
for some arbitrary $\delta>0$. 
Taking  the $\limsup$ in \eqref{inegdelta} as $|x-y|\to 0$ ($x,y\in B(0,R)$) yields
\begin{align*}
\limsup_{|x-y|\to 0}\E^B\big[&\big|e^{-\lambda F(x,\epsilon)}-e^{-\lambda F(y,\epsilon)}\big|\big]\leq \E^B\big[\min\big(2,\lambda  F(x,\delta)+\lambda F(y,\delta)\big)\big].
\end{align*}
It is proved in \cite{GRV} that, almost surely in $X$: 
$$\sup_{x\in B(0,R)}E^{B}[F(x,\delta)]\to 0\quad \text{as }\delta\to 0.$$
Therefore, we can choose $\delta$ arbitrarily close to $0$ to get
\begin{align*}
\limsup_{|x-y|\to 0}\E^B\big[ \big|e^{-\lambda F(x,\epsilon)}-e^{-\lambda F(y,\epsilon)}\big|\big]=0.
\end{align*}
By gathering the above considerations, we have proved that $x\mapsto R_\lambda^{X,\epsilon}f(x)$ is continuous over $D$. Since the family $(R_\lambda^{X,\epsilon}f)_\epsilon$ uniformly converges towards $R^X_\lambda f$ on the compact sets as $\epsilon\to 0$, we deduce that $R^X_\lambda f$ is continuous.
\qed

As a consequence of the above theorem, we can deduce the existence of the Liouville heat kernel:

\begin{theorem}{\bf Liouville heat kernel. }\label{th:heat} For $\gamma\in [0,2[$, the Liouville semigroup $(P^X_\t)_{\t>0}$ is absolutely continuous with respect to the Liouville measure. There exists a family $(\p^X_\t(\cdot,\cdot))_{\t\geq 0}$, called the Liouville heat kernel, of jointly measurable functions such that:
$$\forall f \in B_b(D),\quad P^X_\t f(x)=\int_{D}f(y)\p^X_\t(x,y)\,M(dy)$$ and such that:
\begin{enumerate}
\item (positivity) for all $\t>0$ and for all $x\in D$, for $M(dy)$-almost every $y\in D$, $$\p^X_\t(x,y)\geq 0,$$  
\item (symmetry) for all $\t>0$ and for  every $x,y\in D$: $$\p^X_\t(x,y)=\p^X_\t(y,x),$$
\item (semigroup property) for all $\s,\t\geq 0$, for all $x,y \in D$, $$\p_{\t+\s}^X(x,y)=\int_{D}\p^X_t(x,z)\p^X_{\s}(z,y)\,M(dz).$$  
 \end{enumerate}
\end{theorem}

\vspace{1mm}
\noindent {\it Proof.} 
Since the Liouville semigroup $(P^X_\t)_{\t>0}$ is symmetric with respect to the Liouville measure, we use \cite[Theorem 4.2.4]{fuku} which states that absolute continuity of the resolvent family $R_\lambda^X$ for all $\lambda>0$ is equivalent to absolute continuity of the semigroup. Therefore, it suffices to prove that, almost surely in $X$, 
$$\forall A \text{ Borelian set},\quad M(A)=0\Rightarrow \forall x\in D,\quad R^X_\lambda\ind_A(x)=0.$$
Since the Liouville semigroup is invariant under the Liouville measure, we have for all bounded Borelian set $A$
\begin{equation}\label{resolmeas}
\lambda\int_D R^X_\lambda\ind_A(x)\,M(dx)=M(A).
\end{equation}
Therefore, $M(A)=0$ implies that for $M$-almost every $x\in D$: $R^X_\lambda\ind_A(x)=0$. Since $M$ has full support, we thus have at hand a dense subset $D_A$ of $\R^2$ such that $R^X_\lambda\ind_A(x)=0$ for $x\in D_A$. From Theorem \ref{expr:resol}, the mapping $x\mapsto R^X_\lambda\ind_A(x)$ is continuous. Therefore, it is identically null. Absolute continuity follows.

 \qed

\vspace{2mm}
Now we focus on another aspect of the regularizing properties of the resolvent family (which was already stated in the introduction as theorem \ref{i.holdresolbis}):
\begin{theorem}\label{holdresolbis}
Assume $D=\R^2$ and $\gamma\in [0,2[$. There is an exponent $\alpha\in(0,1)$ (depending only on $\gamma$), such that, almost surely in $X$,  for all $\lambda>0$ the Liouville resolvent is locally $\alpha$-H\"older. More precisely, for each $R$ and $\lambda_0>0$, we can find a random constant $C_{R,\lambda_0}$, which is $\Pb^X$-almost surely finite such that, for all $\lambda\in ]0,\lambda_0]$ and  for all continuous function $f:\R^2\to \R$ vanishing at infinity:
$$\forall x,y\in B(0,R),\quad   |R^X_\lambda f(x)-R^X_\lambda f(y)|\leq \lambda^{-1}C_{R,\lambda_0} \|f\|_\infty|x-y|^\alpha.$$
\end{theorem}

\vspace{1mm}
\noindent {\it Proof.}
Fix $\lambda>0$. Let $f:\R^2\to\R$ be a bounded Borelian function. Let us prove that $x\mapsto R^X_\t f(x)$ is locally H\"older. Without loss of generality, we may assume that $\|f\|_{\infty}\leq 1$. To this purpose, let us work inside a ball centered at $0$ with fixed radius, say $1$. Inside this ball, we consider two different points $x,y$. From these two points, we consider two Brownian motions $B^x$ and $B^y$ coupled in the usual fashion and such that they coincide after some stopping time $\tau^{x,y}$. By applying the strong Markov property, we get:
\begin{align*}
R_\lambda^Xf(x)&=\E^B\Big[\int_0^{\infty}e^{-\lambda F(x,s)}f(B^x_s)\,F(x,ds)\Big]\\
&=\E^B\Big[\int_0^{\tau^{x,y}}e^{-\lambda F(x,s)}f(B^x_s)\,F(x,ds)\Big]+\E^B\Big[\int_{\tau^{x,y}}^{\infty}e^{-\lambda F(x,s)}f(B^x_s)\,F(x,ds)\Big]\\
&\stackrel{def}{=}N_{x,y}(x)+R_\lambda^{X,x,y}f(x).
\end{align*}
 First we focus on $N_{x,y}$:
\begin{align}
|N_{x,y}(x)|\leq &  \E^B\Big[\int_0^{F(x,\tau^{x,y})}e^{-\lambda s}\,ds\Big]\nonumber\\
=& \frac{1}{\lambda}\E^B\Big[1-e^{-\lambda F(x,\tau^{x,y})}\Big]\nonumber\\ 
\leq &\frac{1}{\lambda} \E^B\Big[\min\big(1,\lambda F(x,\tau^{x,y}\big)\Big].\label{est1}
\end{align}
  
Let us now treat the term $R_\lambda^{X,x,y}f$. By the strong Markov property of the Brownian motion, we get:
\begin{align*}
R_\lambda^{X,x,y}f(x)&=\E^B\Big[\int_{\tau^{x,y}}^{\infty}e^{-\lambda F(x,s)}f(B^x_s)\,F(x,ds)\Big]\\
&=\E^B\Big[e^{-\lambda F(x,\tau^{x,y})}R_\lambda^Xf(B^x_{\tau^{x,y}})\Big].
\end{align*}
Therefore we have:
\begin{align*}
|R_\lambda^{X,x,y}f &(x)-R_\lambda^{X,y,x}f(y)|\\=&\Big|\E^B\Big[e^{-\lambda F(x,\tau^{x,y})}R_\lambda^Xf(B^x_{\tau^{x,y}})\Big]-\E^B\Big[e^{-\lambda F(y,\tau^{x,y})}R_\lambda^Xf(B^y_{\tau^{x,y}})\Big]\Big|\\
\leq &\Big|\E^B\Big[e^{-\lambda F(x,\tau^{x,y})}R_\lambda^Xf(B_{\tau^{x,y}}^x)\Big]-\E^B\Big[e^{-\lambda F(x,\tau^{x,y})}R_\lambda^Xf(B_{\tau^{x,y}}^y)\Big]\Big|\\
&+\Big|\E^B\Big[e^{-\lambda F(x,\tau^{x,y})}R_\lambda^Xf(B_{\tau^{x,y}}^y)\Big]-\E^B\Big[e^{-\lambda F(y,\tau^{x,y})}R_\lambda^Xf(B_{\tau^{x,y}}^y)\Big]\Big|\\
=&\Big|\E^B\Big[e^{-\lambda F(x,\tau^{x,y})}R_\lambda^Xf(B_{\tau^{x,y}}^y)\Big]-\E^B\Big[e^{-\lambda F(y,\tau^{x,y})}R_\lambda^Xf(B_{\tau^{x,y}}^y)\Big]\Big|\\
\leq &\lambda^{-1}  \E^B\big[\big|e^{-\lambda F(x,\tau^{x,y})}-e^{-\lambda F(y,\tau^{x,y})}\big|\big].
\end{align*}
In the above inequalities, we have used the facts that $B_{\tau^{x,y}}^x=B_{\tau^{x,y}}^y$ and $\|R_\lambda^Xf\|_{\infty}\leq \lambda^{-1}$. It is readily seen that this  quantity can be estimated by:
\begin{align}
\E^B\big[&\big|e^{-\lambda F(x,\tau^{x,y})}-e^{-\lambda F(y,\tau^{x,y})}\big|\big]\nonumber\\
\leq &\E^B\big[\min\big(1,\lambda | F(x,\tau^{x,y})-F(y,\tau^{x,y})|\big)\big]\nonumber\\
\leq &\E^B\big[\min\big(1,\lambda  F(x,\tau^{x,y})+\lambda F(y,\tau^{x,y})\big)\big]\label{est2}
\end{align}
Therefore, we can take the $q$-th power ($q\geq 1$) and use the Jensen inequality to get
\begin{align}\label{est3}
|R_\lambda^{X,x,y}f  (x)-R_\lambda^{X,y,x}f(y)|^q\leq &C\lambda^{-q}\E^B\big[\min\big(1,\lambda  F(x,\tau^{x,y})+\lambda F(y,\tau^{x,y})\big)^q\big].
\end{align}
Observe that this bound holds for all $\lambda$ and all functions $f$ with $\|f\|_{\infty}\leq 1$. So, let us choose  a countable family $(g_n)_n$ of functions in $C_0(\R^d)$ dense for the topology of uniform convergence over compact sets and set $f_n=g_n/\|g_n\|_\infty$. By gathering \eqref{est1}+\eqref{est3}, we get
$$\E^X\Big[\sup_{\lambda\leq \lambda_0}\sup_n\lambda^q |R_\lambda^{X}f_n  (x)-R_\lambda^{X}f_n(y)|^q\Big]\leq C\E^X\E^B\big[\min\big(1,\lambda_0  F(x,\tau^{x,y})+\lambda_0 F(y,\tau^{x,y})\big)^q\big].$$
We claim:
\begin{lemma}\label{Ftauxy}
For all $x,y\in B(0,1)$ and all $\chi\in]0,\frac{1}{2}[$, $\epsilon>0$, $p\in]0,1[$ and $q\geq 1$ such that $pq>1$, we have
$$\E^X\E^B\Big[\min\big(1,\lambda_0 F(x,\tau^{x,y})\big)^q\Big]\leq C_{\chi,p,q}\Big(\lambda_0^{p q}|x-y|^{(2-\epsilon)\xi(pq) } +|x-y|^{\epsilon q\chi}\Big),$$
for some constant $C_{\chi,p,q}$ which only depends on $\chi,p,q$ and 
$$\forall q\geq 0,\quad \xi(q)=\big(1+\frac{\gamma^2}{4}\big)q-\frac{\gamma^2}{4}q^2.$$
\end{lemma}

 We postpone the proof of this lemma and come back to the proof of Theorem \ref{holdresolbis}. We deduce that for all $x,y\in B(0,1)$, $\chi\in]0,\frac{1}{2}[$, $\epsilon>0$, $p\in]0,1[$ and $q\geq 1$ such that $pq>1$, we have
$$\E^X\Big[\sup_{\lambda\leq \lambda_0}\sup_n\lambda^{q}|R_\lambda^{X}f_n  (x)-R_\lambda^{X}f_n(y)|^q\Big]\leq C_{\chi,p,q,\lambda_0}\Big(|x-y|^{(2-\epsilon)\xi(pq) } +|x-y|^{\epsilon q\chi}\Big),$$
for some constant $C_{\chi,p,q,\lambda_0}$ which only depends on $\chi,p,q,\lambda_0$. Now we fix $\chi\in]0,\frac{1}{2}[$. Then we choose $\delta>0$ such that $1+\delta<\min(2,\frac{4}{\gamma^2})$ (this is possible since $\gamma<2$). Since $\xi(1+\delta)>1$, we can choose $\epsilon>0$ such that $(2-\epsilon)\xi(1+\delta)>2$. Then we choose $q>1$ large enough so as to make $\epsilon \chi q>2$. Then we choose $p\in]0,1[$ such that $pq=1+\delta$. We get
$$\E^X\Big[\sup_{\lambda\leq \lambda_0}\sup_n\lambda^q |R_\lambda^{X}f_n  (x)-R_\lambda^{X}f_n(y)|^q\Big]\leq C_{\chi,p,q,\lambda_0} |x-y|^{\beta }$$
for some $\beta>2$ only depending on $\gamma\in [0,2[$.

From Theorem \ref{th:kolm}, we deduce that for some $\alpha>0$ (only depending on $\gamma$) and some positive $\Pb^X$-almost surely finite random variable $\widetilde{C}$  independent of $n$ and $\lambda\in ]0,\lambda_0]\cap\Q$:
\begin{equation}\label{nique}
\sup_n\sup_{\lambda\in ]0,\lambda_0]\cap\Q}\lambda \big|R_\lambda^{X}g_n  (x)-R_\lambda^{X}g_n(y)\big| \leq  \widetilde{C}\|g_n\|_\infty |x-y|^{\alpha}.\end{equation}
Observe that this relation is then necessarily true for all $\lambda\in ]0,\lambda_0]$ because of the continuity of the resolvent with respect to the parameter $\lambda$.
 Now consider a function $f\in C_0(\R^2)$. There exists a subsequence $(n_k)_k$ such that $\|f-g_{n_k}\|_\infty\to 0$ as $k\to\infty$. In particular $\sup_k\|g_{n_k}\|_\infty<+\infty$ and $\lim_{k\to\infty}\|g_{n_k}\|_\infty=\|f\|_\infty$.  It is plain to deduce from the uniform convergence of $(g_{n_k})_k$ towards $f$ (and therefore the uniform convergence of $R_\lambda^{X}g_n$ towards $R_\lambda^{X}f$) and \eqref{nique} that:
$$\forall x,y\in B(0,R),\quad \lambda\big|R_\lambda^{X}f  (x)-R_\lambda^{X}f(y)\big| \leq \widetilde{C} \|f\|_\infty |x-y|^{\alpha}.$$
The proof is over.\qed

\vspace{1mm}
\noindent {\it Proof of Lemma \ref{Ftauxy}.} Let us consider $R>0$ such that $R|x-y|^2\leq 1$. We have 
\begin{align*}
\E^B\Big[\min&\big(1,\lambda_0 F(x,\tau^{x,y})\big) \Big]\\
=& \E^B\Big[\min\big(1,\lambda_0 F(x,\tau^{x,y})\big) ;\tau^{x,y}\leq R|x-y|^2\Big]\\
&+\E^B\Big[\min\big(1,\lambda_0 F(x,\tau^{x,y})\big);\tau^{x,y}> R|x-y|^2\Big]\\
\leq & \E^B\Big[\min\big(1,\lambda_0 F(x,R|x-y|^2)\big)\Big]+\Pb^B\big(\tau^{x,y}> R|x-y|^2\big)\\
\leq & \E^B\Big[\min\big(1,\lambda^{p }_0 F(x,R|x-y|^2)^{p}\big)\Big]+\Pb^B\big(\tau^{x,y}> R|x-y|^2\big).
\end{align*}
The last inequality results from the fact that $0<p<1$. Therefore, for any $\chi\in]0,\frac{1}{2}[$
\begin{align*}
\E^B\Big[\min&\big(1,\lambda_0 F(x,\tau^{x,y})\big) \Big]\\
\leq & \lambda^{p }_0\E^B\Big[ F(x,R|x-y|^2)^{p } \Big]+\Pb^B\big(\tau^{x,y}> R|x-y|^2\big)\\
\leq & \lambda^{p }_0\E^B\Big[ F(x,R|x-y|^2)^{p } \Big]+C_{\chi}R^{-\chi},
\end{align*}
the last inequality resulting from the fact that the law of the random variable $\tau^{x,y}|x-y|^{-2}$ is stochastically dominated  by a fixed random variable (independent from $x,y$) which possesses moments of order $\chi$ for all $\chi\in ]0,\frac{1}{2}[$. Indeed, if $x=(x_1,x_2)$ and $y=(y_1,y_2)$, then we have the following equality in law $\tau^{x,y} = \max(|y_1-x_1|^{2} \tau_1,|y_2-x_2|^{2}\tau_{2}) $ where $\tau_1$ and $\tau_{2}$ are standard independent stable laws of index $1/2$ (this results from the fact that the hitting time process of a standard $1d$ Brownian motion follows a stable Levy process of index $1/2$: see for instance \cite{Kara}). By taking the $q$-th power and integrating with respect to $\E^X$, we get:
\begin{align*}
\E^X\Big[\E^B\Big[\min&\big(1,\lambda_0 F(x,\tau^{x,y})\big) \Big]^q\Big]\\
\leq & 2^{q-1}\lambda^{p q}_0\E^X\Big[\E^B\Big[ F(x,R|x-y|^2)^{p } \Big]^q\Big]+2^{q-1}C_{\chi}^qR^{-q\chi}\\
\leq & 2^{q-1}\lambda^{p q}_0\E^X \E^B\Big[ F(x,R|x-y|^2)^{pq } \Big]+2^{q-1}C_{\chi}^qR^{-q\chi}.
\end{align*}
Now we take $R=|x-y|^{-\epsilon}$ and use  (see \cite{GRV}):
\begin{proposition}\label{lemmamoments}
If $\gamma^2<4$ and $x\in\R^2$, the mapping $F(x,\cdot)$ possesses moments of order $0\leq q<\min(2,4/\gamma^2)$. 
Furthermore,  if $F$ admits moments of order $q\geq 1$ then, for all $s\in [0,1]$ and $t\in[0,T]$:
$$\E^X\E^B[F(x,[t,t+s])^q]\leq C_qs^{\xi(q)},$$ where $C_q>0$ (independent of $x,T$) and
$$\xi(q)=\big(1+\frac{\gamma^2}{4}\big)q-\frac{\gamma^2}{4}q^2.$$
\end{proposition}
Thus we get:
\begin{align*}
\E^X\Big[\E^B\Big[\min&\big(1,\lambda_0 F(x,\tau^{x,y})\big) \Big]^q\Big]\\
\leq & C_{\chi,p,q}\Big(\lambda^{p q}_0|x-y|^{(2-\epsilon)\xi(pq) } +|x-y|^{\epsilon q\chi}\Big),
\end{align*}  
and we prove the Lemma.\qed

\vspace{1mm}
When $D$ is compact, i.e. when $D=\T^2$ or $\S^2$, we get:
\begin{theorem}\label{holdresolbis2}
Assume $D=\T^2$ or $D=\S^2$ and $\gamma\in [0,2[$. There is an exponent $\alpha\in(0,1)$ (depending only on $\gamma$), such that, almost surely in $X$,  for all $\lambda>0$ the Liouville resolvent is $\alpha$-H\"older. More precisely, for each  $\lambda_0>0$, we can find a random constant $C_{\lambda_0}$, which is $\Pb^X$-almost surely finite such that, for all $\lambda\in ]0,\lambda_0]$ and  for all continuous function $f:D\to \R$, $\forall x,y\in D$:
$$  |R^X_\lambda f(x)-R^X_\lambda f(y)|\leq \lambda^{-1}C_{R} \|f\|_\infty|x-y|^\alpha.$$
\end{theorem}

\vspace{2mm}
\begin{corollary}
For each $\lambda>0$, the resolvent operator $R^X_\lambda:C_0(\R^2)\to C_b(\R^2)$ is compact for the topology of convergence over compact sets.
In the case of the sphere $\S^2$ (or the torus $\T^2$) equipped with a GFF $X$, the resolvent operator $R^X_\lambda:C_b(\S^2)\to C_b(\S^2)$ is compact.
\end{corollary}

\vspace{1mm}
\noindent {\it Proof.} This is just a consequence of Theorems \ref{holdresolbis} or \ref{holdresolbis2}.\qed

\vspace{2mm}
Theorem \ref{expr:resol} has the following consequences on the structure of the resolvent family:
\begin{theorem}{\bf (massive Liouville  Green kernels)}.\label{th:green}
The resolvent family $(R^X_\lambda)_{\lambda>0}$ is absolutely continuous with respect to the Liouville measure. Therefore there exists a family $(\r^X_\lambda(\cdot,\cdot))_\lambda$, called the family of massive Liouville Green kernels, of jointly measurable functions such that:
$$\forall f \in B_b(D),\quad R^X_\lambda f(x)=\int_{D}f(y)\r^X_\lambda(x,y)\,M(dy)$$ and such that:\\
1) (lower semi-continuity) For all $x \in D$, the function $y \mapsto \r^X_\lambda(x,y)$ is lower semi-continuous. $$$$
2) (strict-positivity) for all $\lambda>0$ and for all compact set $K$, $$\inf_{x,y \in K} \r^X_\lambda(x,y)>0,$$   
3) (symmetry) for all $\lambda>0$ and for  every $x,y\in D$: $$\r^X_\lambda(x,y)=\r^X_\lambda(y,x),$$
4) (resolvent identity) for all $\lambda,\mu>0$, for all $x,y \in D$, $$\r^X_\mu(x,y)-\r^X_\lambda(x,y)=(\lambda-\mu)\int_{D}\r^X_\lambda(x,z)\r^X_{\mu}(z,y)\,M(dz).$$  
5) ($\lambda$-excessive) for every $y$, 
$$e^{-\lambda t}P_t^X(\r_\lambda(\cdot,y))(x)\leq \r_\lambda(x,y)$$ for $M$-almost every $x$ and for all $t>0$.
\end{theorem}
 
\vspace{1mm}
\noindent {\it Proof.} Since the Liouville semigroup is absolutely continuous with respect to the Liouville measure (see theorem \ref{th:heat}), we can apply \cite[Lemma 4.2.4]{fuku} which proves the existence of the massive Liouville Green kernels with the items 3), 4) and 5). Next, we check item 1). In fact, in \cite[Lemma 4.2.4]{fuku}, the resolvent is constructed such that for all $x,y \in D$
\begin{equation*}
 \r_\lambda(x,y)= \underset{t \to 0}{\lim} \: e^{-\lambda t}P_t^X(\r_\lambda(\cdot,y))(x)=  \underset{t \to 0}{\lim}\:  R^X_\lambda  ( e^{-\lambda t} \p^X_\t(x, \cdot)  ) (y)
\end{equation*} 
where the above limit is non-decreasing as $t$ goes to $0$. Therefore, by theorem \ref{expr:resol}, the resolvent density is such that, for all $x \in D$, the function $y \mapsto \r^X_\lambda(x,y)$ is lower semi-continuous.  

Finally, we check item 2). Let $\lambda>0$ and some compact set $K$ be fixed. First, consider a Borel set $A$ and $x\in D$ such that $R_\lambda^X\ind_A(x)=0$. Then $\Pb^B_x$ a.s., we have $$\int_0^{\infty}e^{-\lambda F(x,t)}\ind_A(B_t^x)\,F(x,dt)=0.$$
We deduce, by using the Markov property:
\begin{align*}
0=& \E^B_x\Big[\int_1^{\infty}e^{-\lambda F(x,t)-F(x,1)}\ind_A(B_t^x)\,F(x,dt)\Big]\\
=&\E^B_x\Big[R_\lambda^X\ind_A(B_1^x)\Big].
\end{align*}
Since the transition probabilities of the standard Brownian motion are strictly positive, we deduce that $\Pb^X$ almost surely the mapping $x\mapsto R_\lambda^X\ind_A(x)$ vanishes over a set with full Lebesgue measure. Furthermore, it is continuous by Theorem \ref{expr:resol}. Thus we have $R_\lambda^X\ind_A=0$ identically. Finally, we get:
$$M(A)= \lambda \int_D R_\lambda^X\ind_A(x)\,M(dx)=0.$$
Since this is true for all Borel sets $A$, we deduce that $\Pb^X$-almost surely the resolvent density $\r^X_\lambda(x,\cdot)$ is positive $M$ almost everywhere.

Now, suppose that one can find a sequence $(x_y,y_n)_{n \geq 1}$ in $K$ such that $\r^X_\lambda(x_n,y_n)$ converges to $0$ as $n$ goes to infinity. By compactness of $K$, one can assume that the sequence $(x_y,y_n)_{n \geq 1}$ converges to $(x,y)$. Now, we have by item 4) that    
\begin{equation*}             
\r^X_\lambda(x_n,y_n)    \geq \lambda \int_{D}\r^X_\lambda(x_n,z)\r^X_{2 \lambda}(z,y_n)\,M(dz).
\end{equation*}
Taking the limit above as $n$ goes to infinity leads by Fatou's lemma and item 1) to 
\begin{equation*}             
\underset{n \to \infty}{\underline{\lim}}\r^X_\lambda(x_n,y_n)    \geq \int_{D} \underset{n \to \infty}{\underline{\lim}} \r^X_\lambda(x_n,z) \underset{n \to \infty}{\underline{\lim}}\r^X_{2 \lambda}(z,y_n)\,M(dz) \geq \int_{D}  \r^X_\lambda(x,z)\r^X_{2 \lambda}(z,y)\,M(dz).
\end{equation*}
Therefore $\int_{D}  \r^X_\lambda(x,z)\r^X_{2\lambda}(z,y)\,M(dz)=0$ which contradicts the fact that $\r^X_\lambda(x,\cdot)$ and $\r^X_\lambda(\cdot,y)$ are positive $M$ almost everywhere. Hence, we get $\inf_{x,y \in K} \r^X_\lambda(x,y)>0$.\qed

\vspace{0.5cm}

We end this section by collecting some consequences of the above analysis on the behavior of $\LB_\t$ with respect to the Liouville measure $M$.

\subsection{Recurrence and ergodicity}

As prescribed in \cite[section 1.5]{fuku}, let us define the Green function for $f\in L^1(D,M)$  by (if exists)
$$ Gf(x)=\lim_{t\to\infty} \int_0^tP^X_rf(x)\,dr.$$
We further denote $g_D$ the standard Green kernel on $D$. 

Following \cite{fuku}, we say that the semi-group $(P^X_t)_t$, which is symmetric w.r.t. the measure $M$,   is   {\it irreducible} if  any $P^X_t$-invariant set $B$ satisfies $M(B)=0$ or $M(B^c)=0$. We say that $(P_t^X)$ is recurrent if, for any $f\in L^1_+(D,M)$, we have $Gf(x)=0$ or $Gf(x)=+\infty$ $M$-almost surely.

\begin{theorem}{\bf (Liouville Green function)} The Liouville semi-group is irreducible and recurrent. If $D$ is the torus $\T^2$ or the sphere $\S^2$, the Liouville Green function, denoted by $G^X_D$, is given by
$$G^X_Df(x)=\int_Dg_D(x,y)f(y)\,M(dy)$$ for all functions $f\in L^2(D,M)$ such that
$$\int_Df(y)\,M(dy)=0.$$
\end{theorem}

\vspace{1mm}
\noindent {\it Proof.} Irreducibility is a straightforward consequence of Theorem \ref{th:green}. 

Let us establish recurrence. We carry out the proof in the case of $D=\R^2$. The reader can easily adapt the proof to $D=\S^2,\T^2$. We first observe that
\begin{align*}
Gf(x)=&\lim_{t\to\infty} \int_0^tP^X_rf(x)\,dr\\
=&\lim_{t\to\infty}\E^B_x\Big[ \int_0^tf(\LB^x_\r)\,d\r\Big]\\
=&\lim_{t\to\infty}\E^B_x\Big[ \int_0^{F^{-1}(x,t)}f(B_r^x)\,F(x,dr)\Big]\\
=&\lim_{t\to\infty}\E^B_x\Big[ \int_0^{t}f(B_r^x)\,F(x,dr)\Big].
\end{align*}
We have used the fact that $F(x,t)$ almost surely converges towards $+\infty$ as $t\to\infty$ (see \cite{GRV}). Now observe that
\begin{align*}
\E^B\Big[ \int_0^{t}f(B_r^x)\,F(x,dr)\Big]=& \int_{\R^2}\Big(\int_0^t\frac{e^{-\frac{|x-u|^2}{2r}}}{2 \pi r}\,dr\Big) f(u)\,M(du).
\end{align*}
Let us assume that $f$ has compact support. The above quantity diverges as $t\to \infty$ like $C\ln t$. To see this, we first compensate the divergence at infinity as follows:
\begin{align*}
\E^B_x\Big[ \int_0^{t}f(B_r^x)\,F(x,dr)\Big]=& \int_{\R^2}\Big(\int_0^t\frac{e^{-\frac{|x-u|^2}{2r}}-e^{-\frac{1}{2r}}}{2 \pi r}\,dr\Big) f(u)\,M(du)\\
&+\Big(\int_0^t\frac{e^{-\frac{1}{2r}}}{2 \pi r}\,dr\Big)\int_{\R^2}f(u)\,M(du).
\end{align*}
By passing to the limit as $t\to\infty$, we get
\begin{align}
\lim_{t\to\infty}\int_{\R^2}&\Big(\int_0^t\frac{e^{-\frac{|x-u|^2}{2r}}-e^{-\frac{1}{2r}}}{2 \pi r}\,dr\Big) f(u)\,M(du)\nonumber\\
&=\int_{\R^2}\Big(\int_0^{\infty}\frac{e^{-\frac{|x-u|^2}{2r}}-e^{-\frac{1}{2r}}}{2 \pi r}\,dr\Big) f(u)\,M(du).\label{eq:green}
\end{align}
In fact, the above computations need some further explanations. The Green kernel appears through the relation $$\Big(\int_0^{\infty}\frac{e^{-\frac{|x-u|^2}{2r}}-e^{-\frac{1}{2r}}}{2 \pi r}\,dr\Big)=\frac{1}{\pi}\ln\frac{1}{|u-x|}.$$ To prove that this log term does not affect the convergence of the integral \eqref{eq:green}, we have to use Corollary 2.21 in \cite{GRV}.  

Obviously, we are thus left with two options. Either
$$\int_{\R^2}f(u)\,M(du)=0$$ and $f=0$ $M$ almost surely, which entails $Gf=0$ $M$-almost surely, or
$$\int_{\R^2}f(u)\,M(du)>0$$ leading to $Gf(x)=+\infty$. Put in other words, the Liouville semi-group is recurrent. 

Now let us suppose that $D$ is the torus $\T^2$ or the sphere $\S^2$. The expression for the Green function will rely heavily on \cite[section 4.8]{fuku} and more specifically on Theorem 4.8.1 of \cite{fuku}. In order to apply  Theorem 4.8.1, one must check the condition (4.8.3) of \cite{fuku}, i.e. that there exists some constant $c>0$ such that $\r^X_1(x,y) >c$ for all $x,y \in D$. This is exactly the content of item 2) in Theorem \ref{th:green}. Now Theorem 4.8.1 of \cite{fuku} implies that for all $f\in L^2(D,M)$ such that $\int_Df(y)\,M(dy)=0$, $R^X_\lambda(f)$ converges in $L^2(D,M)$ to $\int_Dg_D(x,y)f(y)\,M(dy)$ as $\lambda$ goes to $0$. Theorem 4.8.1 of \cite{fuku} also implies a Poincar\'e type inequality: as a consequence, $Gf$ also converges in $L^2(D,M)$ to $\int_Dg_D(x,y)f(y)\,M(dy)$.

\qed

\begin{remark}
It may be worth saying here that the above proof works in dimension $2$ (or $1$) only. It is based on the recurrence of the standard $2$-dimensional Brownian motion. This is one point were techniques related to the Liouville Brownian motion differ according to the dimension.

Let us also point out that  the integral formula for the Liouville Green function should be convenient to study its regularizing property. For instance, it is almost obvious to see that $Gf$ is continuous when $f$ is bounded.
\end{remark}

We can now apply \cite[Theorem 4.7.3]{fuku} to get
\begin{theorem}
Let us denote by $\Pb_M$ the law of the Liouville Brownian motion with initial distribution $M$.  Let $f\in L^1(D,M)$ be a Borel measurable function.
\begin{enumerate}
\item It holds  $\Pb_M$-almost surely that
$$ \lim_{t\to\infty}\frac{1}{t}\int_0^tf(\LB_\r)\,d\r=\frac{1}{M(D)}\int_Df(x)\,M(dx).$$
\item Assume further that $f$ is locally uniformly bounded. Then, $\Pb_x$-almost surely 
$$ \lim_{t\to\infty}\frac{1}{t}\int_0^tf(\LB_\r)\,d\r=\frac{1}{M(D)}\int_Df(x)\,M(dx).$$
\end{enumerate}
\end{theorem}
\begin{remark}
Observe that the above theorem entails uniqueness of the invariant probability measure for the Liouville Brownian motion in the case of the sphere or the torus. With additional efforts, one could also establish in this way uniqueness in the case of the whole plane up to multiplicative constants.
\end{remark}

\subsection{The Liouville Brownian motion spends most of his time in the thick points of $X$}

In this section, it is convenient to make explicit the dependence of $\gamma$ of the Liouville measure, i.e. we write $M_\gamma$ instead of $M$. Following Kahane \cite{cf:Kah}, let us introduce the $\gamma$-thick points of $X$ (see also \cite{BMI} in the one dimensional case and \cite{hmp} in the case of circle averages of the GFF):
 \begin{equation*}
 K_\gamma= \lbrace x \in D; \: \underset{n \to \infty}{\lim} \: \frac{  X_n(x)}{\ln  c_{n+1}} = \gamma  \rbrace \,, 
 \end{equation*}
 where the series $(c_n)_{n\geq 1}$ was introduced in \cite{GRV}, i.e. $$\E[X_n(x) X_n(y)]=\int_{1}^{c_{n+1}}\frac{k_m(u(x-y))}{u}\,du.$$
 
 It is a well known fact that the Borel set $K_\gamma$ gives full mass to the measure $M_\gamma$, i.e. $M_\gamma(K^c_\gamma)=0$: this was proved in Kahane's seminal work \cite{cf:Kah} using the so-called Peyri\`ere measure. The sets $ K_\gamma$ appear also frequently in the general context of multifractal formalism for multifractal measures (log-Poisson, discrete cascades, etc...). The terminology "thick points" is not due to Kahane but appears in \cite{hmp} for example. This implies that, if $\gamma, \gamma' \in  [0,2[$ are such that $\gamma \not = \gamma'$, the measures $M_\gamma$ and $M_{\gamma'}$ are singular with respect to each other. In particular, the measure $M_\gamma$ for $\gamma\in ]0,2[$ is singular with respect to the Lebesgue measure (which corresponds to $\gamma=0$). Let us also stress that $K_\gamma$ should be distinguished from the support of $M_\gamma$ which is $D$.

 As a consequence of Theorem \ref{th:green}, we obtain the following result where $\lambda$ is the Lebesgue measure: 
\begin{corollary}
For $\gamma\in [0,2[$, the Liouville Brownian motion spends Lebesgue-almost all its time in the $\gamma$-thick points of $X$  for all starting points $x$: 
 \begin{equation*}
\text{a.s. in }X, \forall x \in D, \text{ a.s. under }\Pb_x^{B},\quad \lambda \lbrace \t \geq 0; \:  \LB^x_\t \in    K^c_\gamma  \rbrace =0.
\end{equation*}
\end{corollary}

\begin{remark}
\noindent
\bi
\item A weaker form of this result is proved in \cite{berest}: the author proves that  $$\lambda \lbrace \t \geq 0; \:  \LB^x_\t \in     K^c_\gamma  \rbrace =0$$

for  one fixed starting point.
\item This result may also be recovered from the invariance of the measure $M$ proved in \cite{GRV}, but in that case only for M(dx)-almost all starting points (which is thus also slightly weaker than our Corollary). 
\ei
\end{remark}

If one now relies on Theorem \ref{th:heat} instead, one obtains the following Corollary (note that it different from the above one, not stronger, nor weaker). 

\begin{corollary}
For $\gamma\in [0,2[$, almost surely in $X$, for all $\t>0$
$$\Pb_x^{B} \text{a.s.},\quad \LB^x_\t\in K_\gamma.$$
\end{corollary}

Observe that the above corollary was already known when the initial law of the Liouville Brownian motion is the Liouville measure \cite{GRV}. Replacing the starting law by the Dirac mass at $x$ (a.s. in $X$ for all $x$) is a much stronger statement.

\section{Degenerescence of the intrinsic metric associated to the Liouville Dirichlet form}

\subsection{Background on the geometric theory of Dirichlet forms and extension of Riemannian geometry}

As a strongly local regular Dirichlet form, $(\Sigma,\mathcal{F})$ can be written as 
\begin{equation}
\Sigma(f,g)=\int_{\R^2}\,d \Gamma(f,g)
\end{equation}
where $\Gamma$ is a positive semidefinite, symmetric bilinear form on $\mathcal{F}$ with values in the signed Radon measures on $\R^2$ (the so-called energy measure). Denoting by $P_t(x,dy)$ the transition probabilities of the semi-group,  the energy measure can be defined by the formula 
\begin{align*}
\int_{D}\phi \,d\Gamma(f,f)&=\Sigma(f,\phi f)-\frac{1}{2}\Sigma(f^2,\phi)\\
&=\lim_{t\to 0}\frac{1}{2t}\int_{D}\int_{\R^2}\phi(x)(f(x)-f(y))^2 P_t(x,dy)M(dx)
\end{align*}
for every $f\in\mathcal{F}\cap L^\infty(D,M)$ and every $\phi\in\mathcal{F}\cap C_c(D)$. The energy measure is local, satisifes the Leibniz rule as well as the chain rule \cite{fuku}. Let us denote by $\mathcal{F}_{loc}=\{f\in L^2_{loc}(D,M);\Gamma(f,f)\text{ is a Radon measure}\}$. 

The energy measure defines in an intrinsic way a  distance in the wide sense $d_X$ on $D$ by
\begin{equation}\label{defsupremumdist}
d_X(x,y)=\sup\{f(x)-f(y);f\in \mathcal{F}_{loc}\cap C(D),\Gamma(f,f)\leq M \}
\end{equation}
called intrinsic metric \cite{BM1,BM2,DA,VSC}. The condition $\Gamma(f,f)\leq M$ means that the energy measure $\Gamma(f,f) $ is absolutely continuous w.r.t to $M$ with Radon-Nikodym derivative $\frac{d}{dM}\Gamma(f,f)\leq 1$. In general, $d_X$ may be degenerate $d_X(x,y)=0$ or $d_X(x,y)=+\infty$ for some $x\not = y$.

\subsection{Why it vanishes in the setting of Liouville quantum gravity}

Here we provide a rigorous proof in the next subsection followed by  a more heuristical explanation by considering the intrinsic metric associated to the $n$-regularized Dirichlet forms $(\Sigma^n,\mathcal{F}^n)$ obtained by using $X_n$.

\subsubsection{A proof that the intrinsic metric vanishes}

\begin{proposition}
For $\gamma\in [0,2[$, almost surely in $X$, the distance in the wide sense $d_X$ reduces to $0$ for all points $x,y\in D$.
\end{proposition}  

\vspace{1mm}
\noindent {\it Proof.} For $f\in  \mathcal{F}_{loc}\cap C(D)$, the energy measure is characterized by  (see \cite[(3.2.14) and Th. 6.2.1]{fuku}):
$$\forall \phi \in C_c(D), \quad \int_{D}\phi \,d\Gamma(f,f)=2\Sigma (f \phi,f)-\Sigma(f^2,\phi).$$
It is worth mentioning here that the above formula implicitly implies that the energy measure of $f\in L^2(D,M)$ does not depend on the choice of the element $\tilde{f}\in H^1_{loc}(D,dx)$ such that $f=\tilde{f}$ $M$-almost everywhere.

Routine computations on differentiation then entail that
$$\forall \phi \in C_c(D), \quad \int_{D}\phi \,d\Gamma(f,f)=\int_{D}\phi(x) |\nabla f(x)|^2\,dx.$$
Since $M$ and the Lebesgue measure are singular with respect to each other (see above), the condition $\Gamma(f,f)\leq M$ entails that $\nabla f=0$. Therefore, the set $\{f\in \mathcal{F}_{loc}\cap C(D),\Gamma(f,f)\leq M \}$ is the set of constant functions on $D$. Therefore, if $x,y$ are two points in $D$, the supremum in \eqref{defsupremumdist} is $0$. \qed

\subsubsection{A heuristical justification by looking at the $n$-regularized forms}

It is tempting to write in a loose sense that
\begin{equation}\label{}
(\Sigma^n,\mathcal{F}^n) \to (\Sigma, \mathcal{F})\,.
\end{equation}
Now, it is easy to check that the intrinsic metric $d_n$ associated to the Dirichlet form 
\begin{equation}\label{}
\Sigma^n (f,f):= \frac 1 2 \int_D |\nabla f(x)|^2 dx \,,
\end{equation}
with domain $$\mathcal{F}=\Big\{f\in L^2(D,M_n); \nabla f\in L^2(D,dx)\Big\},$$
is exactly the Riemannian distance with metric tensor given by 
\begin{equation}\label{}
g_n(x)   = e^{\gamma X_n(x) - \frac {\gamma^2} 2 \Eb{X_n^2}} dx^2.
\end{equation}
We have the following result, which is in some sense folklore within the community but to our knowledge is not written down anywhere:
\begin{proposition} 
The couple $(D, d_n)$ converges towards the trivial distance, meaning that for all $x,y\in D$, a.s. in $X$:
$$d_n(x,y)\leq C_{x,y}e^{-\frac{\gamma^2}{8}\E[X_n^2(s)]}$$ for some random constant $C_{x,y}>0$.
\end{proposition}

\vspace{1mm} 
\noindent {\it Proof.} By definition, 
$$d_n(x,y)=\inf\{\int_0^1e^{\frac{\gamma}{2}X_n(\sigma_t)-\frac{\gamma^2}{4}\E[X_n^2(\sigma_t)]}|\dot{\sigma}_t|\,dt;\sigma\text{ rectifiable from }x\text{ to }y\}.$$
Obviously, this distance is bounded from above by the weight of the segment joining $x$ to $y$, i.e.
\begin{align*}
d_n(x,y)\leq &\int_{[x,y]}e^{\frac{\gamma}{2}X_n(s)-\frac{\gamma^2}{4}\E[X_n^2(s)]}\,ds\\
=&e^{-\frac{\gamma^2}{8}\E[X_n^2(s)]}\int_{[x,y]}e^{\frac{\gamma}{2}X_n(s)-\frac{\gamma^2}{8}\E[X_n^2(s)]}\,ds
\end{align*}
 where $ds$ stands for the standard arc length on $D$. The arc length restricted to the segment $[x,y]$ is a Radon measure in the class $R_1^+$ of \cite{cf:Kah}. Therefore, for $\gamma\in�[0,2[ $, the limit 
 $$C_{x,y}=\lim_{n\to \infty}\int_{[x,y]}e^{\frac{\gamma}{2}X_n(s)-\frac{\gamma^2}{8}\E[X_n^2(s)]}\,ds$$ exists and is non trivial.\qed

\subsubsection*{Conclusion}

Roughly speaking, one may say that  the geometric aspect of Dirichlet forms at the level of constructing a distance is not as powerful as one might hope looking at its degree of generality. If the machinery seems to be efficient when the underlying space is not too far from a smooth Riemannian geometry, it does not overcome the issue of renormalization. The reader may object that it is not clear that the Liouville distance exists and this could be an explanation to the fact that the intrinsic metric of Dirichlet forms vanishes: we stress that this objection is not relevant   since it does not even work in dimension $1$ though the distance is perfectly explicit and non trivial.


\appendix

\section{Index of notations}\label{index}

\bi
\item $X$: Gaussian Free Field,
\item $M$ (or $M_\gamma$): Liouville measure, 
\item $(B^x_t)_t$: a standard Brownian motion starting from $x$; the corresponding probability measure (expectation) will be denoted $\Pb^B_x$ ($\E^B_x$).
\item  $(\Sigma, \mathcal{F})$: Dirichlet-form,  
\item $(R^X_\lambda)_{\lambda\geq 0}$: Liouville resolvent operator,
\item $B_b(D)$: space of bounded measurable functions on $D$,
\item $C(D)$: space of continuous functions on $D$,
\item $C_b(D)$: space of bounded continuous functions on $D$,
\item $C_0(D)$: space of continuous functions on $D$ vanishing at infinity,
\item $C_c(D)$: space of continuous  functions on $D$ with compact support,
\item $L^p(D,\mu)$: Borel measurable functions on $D$ with $\mu$-integrable $p$-th power, 
\item $H^1(D,dx)$: standard Sobolev space,
\item $H^1_{loc}(D,dx)$: functions which are locally in $H^1(D,dx)$. 

\ei

\section{Reinforced Kolmogorov's continuity criterion}\label{app:kolm}

 In this section, we prove the following result:
 
 \begin{theorem}\label{th:kolm}
 Assume that $(f_n)_n$ is a sequence of random functions defined on the same probability space $(\Omega,\mathcal{F},\Pb)$ such that for some $q,\beta>0$ and for all $x,y\in B(0,R)\subset \R^d$:
 $$\E\big[\sup_n|f_n(x)-f_n(y)|^q\big]\leq C|x-y|^{d+\beta}.$$
For all $\alpha\in]0,\frac{\beta}{q}[$, we can find a modification of $f_n$ for each $n$ (still denoted by $f_n$) and a random constant $\widetilde{C}$, which is $\Pb$-almost surely finite such that:
$$\forall n,\forall x,y\in B(0,R),\quad |f_n(x)-f_n(y)|\leq \widetilde{C}|x-y|^\alpha.$$
 \end{theorem}
 
 \vspace{2mm}
 \noindent {\it Proof.} For simplicity, we carry out the proof in dimension $1$ and we assume that $x,y$ belong to the set $[0,1]$. Let us consider $\alpha \in ]0, \beta/q[$. We get:
\begin{align*}
\Pb\Big( \max_{k=1\dots 2^N}&\sup_n\big|f_n(\frac{ k}{2^N})-f_n(\frac{ k-1}{2^N})\big|>2^{-N\alpha}\Big) \\
& =   \Pb\Big( \bigcup_{
k=1}^{2^N}\Big\{\sup_n\big|f_n(\frac{ k}{2^N})-f_n(\frac{ k-1}{2^N})\big|>2^{-N\alpha}\Big\}\Big) \\ 
& \leq   \sum_{k=1}^{2^N}\Pb\Big( \Big\{\sup_n\big|f_n(\frac{ k}{2^N})-f_n(\frac{ k-1}{2^N})\big|^q >2^{-Nq \alpha }\Big\}\Big)
\\ & \leq  2^{Nq \alpha } \sum_{k=1}^{2^N}\E\Big(\sup_n\big|f_n(\frac{ k}{2^N})-f_n(\frac{ k-1}{2^N})\big|^q
\Big)   \\
 & \leq  2^{Nq \alpha
} \sum_{k=1}^{2^N}C2^{-N(1+\beta) }  \\ 
& =   C  2^{Nq \alpha  -N\beta  }  .
\end{align*} 
Since $\alpha \in ]0,  \beta/q[$, we have $\sum_{N=1}^{\infty }2^{Nq \alpha  -N\beta  } <\infty $. Borel-Cantelli's lemma yields
$$\Pb\Big( \limsup_N \Big\{\max_{k=1\dots 2^N} \sup_n\big| f_n(\frac{ k}{2^N})-f_n(\frac{ k-1}{2^N})\big|>2^{-N\alpha}\Big\}\Big)=0.$$ 
Put in other words, there exists a measurable set $A\in {\cal
F}$ such that $\Pb(A)=1$ and $\forall \omega \in A$, $\exists N_\omega \in \N$, $\forall N\geq N_\omega$, $$\max_{k=1\dots 2^N}\sup_n\big| f_n(\frac{ k}{2^N})-f_n(\frac{ k-1}{2^N})\big|\leq 2^{-N\alpha}.$$
Let us denote by $D_m=\left\{\frac{ k}{2^m},0\leq k\leq
2^m\right\}$ the set of dyadic numbers of order $m$. Let $m,p\in \N $ tels que $m>p\geq N_\omega $, and consider $s,t\in D_m$ such that $s<t$ and $|t-s|\leq 2^{-p}$. 
Then  $s=\frac{ k}{2^m}$ and we can find $a_1,\dots ,a_{m-n}\in \{0,1\}$ such that $t=\frac{ k}{2^m}+\frac{a_1}{2^{p+1}}+\dots+ \frac{ a_{m-p}}{2^m}$. We obtain for $\omega\in A$:
\begin{align*}
\sup_n&|f_n(t,\omega)-f_n(s,\omega )| \\
& =  \sup_n |f_n(\frac{ k}{2^m}+\frac{a_1}{2^{p+1}}+\dots+ \frac{ a_{m-p}}{2^m},\omega)-f_n(\frac{
k}{2^m},\omega )| \\
& \leq   \sum_{j=1}^{m-p}  \sup_n |f_n(\frac{ k}{2^m}+\frac{
a_1}{2^{p+1}}+\dots+ \frac{ a_{j}}{2^{p+j}},\omega )-f_n(\frac{ k}{2^m}+\frac{
a_1}{2^{p+1}}+\dots+ \frac{ a_{j-1}}{2^{p+j-1}},\omega )| \\ & \leq  
\sum_{j=1}^{m-p}  2^{-(p+j)\alpha } .
\end{align*}
Let us now consider $s,t\in D=\bigcup_m D_m$ such that  $|s-t|\leq 2^{-N_\omega}$. Let $p\in \N $ such that $|s-t|\leq 2^{-p}$ and $|s-t|> 2^{-p-1}$. Let $m>p$ such that $s,t\in D_m$. From the previous computations, we get:
\begin{align*}
\sup_n|f_n(t,\omega )-f_n(s,\omega )| & \leq   \sum_{j=p+1}^{m}\frac{
1}{2^{j\alpha }} \\ & \leq  \frac{ 2^\alpha}{2^\alpha
-1}|t-s|^\alpha.
\end{align*}
For each $n$, the mapping $t\mapsto f_n(t,\omega )$ is therefore $\alpha$-H\"older on $D\cap [0,1]$ so that it can be extended to the whole $[0,1]$ while remaining  $\alpha$-H\"older with the same H\"older constant.  Since $f_n$ is continuous in probability for each $n$, this extension is a modification of $f_n$ for all $n$.\qed





\ACKNO{The authors wish to thank G. Miermont for very enlightening discussions
  about the degenerescence of the {\it intrinsic metric} as well as the
  Liouville Green function. We would also like to thank M. De La Salle, A.
  Guillin, M. Hairer, J. Mattingly for fruitful discussions.}


\end{document}